\RequirePackage{fix-cm}
\documentclass[11pt]{article}
\usepackage{graphicx}
\usepackage{amsmath,amsfonts,amssymb,mathtools,dsfont}
\usepackage{xcolor}
\usepackage{enumitem}
\newtheorem{prop}{\textbf{Theorem}}

\begin{document}

\title{A review on asymptotic inference in stochastic differential equations with mixed-effects.}


\author{Maud Delattre\\ Universit\'e Paris-Saclay, INRAE, MaIAGE, 78350, Jouy-en-Josas, France}


\date{}

\maketitle

\begin{abstract}
This paper is a survey of recent contributions on estimation in stochastic differential equations with mixed-effects. These models involve $N$ stochastic differential equations with common drift and diffusion functions but random parameters that allow for differences between processes. The main objective is to estimate the distribution of the random effects and possibly other fixed parameters that are common to the $N$ processes. While many algorithms have been proposed, the theoretical aspects related to estimation have been little studied. This review article focuses only on theoretical inference for stochastic differential equations with mixed-effects. It has so far only been considered in some very specific classes of mixed-effect diffusion models, observed without measurement error, where explicit estimators can be defined. Within this framework, the asymptotic properties of several estimators, either parametric or nonparametric, are discussed. Different schemes of observations are considered  according to the approach, associating a large number of individuals with, in most cases, high-frequency observations of the trajectories. 
\end{abstract}



\section{Introduction}
\label{intro}

The analysis of dynamical phenomena is common in many fields of application such as agronomy with the study of animal or plant growth curves, pharmacology with the study of drug kinetics or even econometrics with the monitoring of productions over time. The available data then often consist of repeated measurements of a continuous-time process among a population of individuals. In biology and in the medical field, data of this kind are referred to as longitudinal data (see \cite{liu2015} for instance), whereas in other disciplines such as econometrics, it is called panel data (see \cite{sul2019}). Whatever the field of study, mixed-effects models are mainly used to analyze simultaneously repeated observations from several individuals. They are also called population models in the sense that they are intended both to describe the variability of a given phenomenon within a population and its specificities for each member of the population individually. Nonlinear mixed-effects models defined by ordinary differential equations are used extensively to describe dynamical processes based on longitudinal data (see \cite{pinheiro2009} for a few examples). In recent years, a number of works have promoted the use of stochastic differential equations (SDEs) with mixed effects as a more realistic alternative to classical nonlinear mixed effects models. In \cite{donnet2010}, Donnet et al. have proposed a SDE version of a Gompertz curve to model the body weights of chicken. In \cite{picchini2019}, Picchini and Forman studied a geometric Brownian motion with random effects in the context of tumor growth modeling. SDEs with mixed-effects have also become increasingly popular in pharmacokinetics and in neuroscience (see \cite{dion2016a} and \cite{donnet2013} for a few examples). This has motivated the development of many algorithms for the estimation of these models. To cite a few examples: Donnet and Samson (2008, 2014) and Delattre and Lavielle (2013) implemented several variants of the SAEM algorithm (\cite{donnet2008}, \cite{donnet2014} and \cite{delattre2013k} respectively) to achieve maximum likelihood estimation, particle algorithms have been proposed in \cite{donnet2014} and \cite{botha2019}, and Bayesian inference was also considered, as in the works of Donnet et al. (2010, \cite{donnet2010}) and Picchini and Forman (2019, \cite{picchini2019}) cited just above. We refer the reader to \cite{picchiniWebpage} for an exhaustive bibliography on inference in SDEs with mixed-effects. Although many methodological developments have accompanied this growing interest in SDEs with mixed-effects, inference for these models has received very little theoretical attention. The main goal of this paper is to present a state of the art of the existing inference methods and to discuss the properties that can be expected for them. Therefore, the focus is here on asymptotic properties of the estimators rather than on computational aspects. 

Section 2 further introduces notations through a general presentation of SDEs with mixed-effects as an extension of classical nonlinear mixed-effects models. The asymptotic properties of the maximum likelihood estimator in nonlinear mixed-effects models are recalled to better position the theoretical inference results obtained in SDEs with mixed-effects. In Section 3, the exact likelihood is given and the technical difficulties that arise with maximum likelihood estimation are briefly discussed in both situations where the individual processes are continuously and discretely observed. Section 4 deals with estimators built from the continuous-time likelihood of the mixed-effects processes. This is the first estimation method to have been studied theoretically. In Section 5, contrasts based on the Euler scheme of the processes are discussed. Section 6 and Section 7 are respectively about plug-in and non-parametric estimation techniques. These latter two approaches rely on preliminary estimations of the random effects. The paper finishes with some discussions.  

\section{Preliminary notes}

\subsection{From nonlinear mixed-effects models to mixed-effects diffusion models}

Let $N$ denote the number of subjects and let $Y_i=(Y_{i1},\ldots,Y_{in_i})'$ be the vector of observations for subject $i$ where $Y_{ij}$ denotes the $j$th measure observed under time condition $t_{ij}$ in individual $i$. For simplicity, we assume in this paper that the number of individual observations is the same for each subject, \textit{i.e.} $n_i = n$ $\forall i \in 1,\ldots,N$. When one is interested in modeling several dynamics over time, there are two main objectives: \textit{i)} provide the best possible description of how the phenomenon under study evolves over time and \textit{ii)} best describe the differences between the individuals within the population of interest. Taking simultaneously into account these two sources of variability leads generally to better statistical results than analyzing each individual dynamics separately. As mentioned above, nonlinear mixed-effects models are classically used for this purpose. Nonlinear mixed effects models can be seen as extensions of standard nonlinear regression models involving the same regression function $f$ for the $N$ individual series of observations, but different parameter values for each individual. More precisely, these are two-level hierarchical models where 
\begin{enumerate}
\item at first level, the $j$th observation of the $i$th subject is modeled as
\begin{equation}
Y_{ij} = f(\psi_i,\beta,t_{ij}) +\varepsilon_{ij} \, , \, i=1,\ldots,N \, , \, j=1,\ldots,n \, ,
\label{eq:nonlinear}
\end{equation}
where the $\varepsilon_i=(\varepsilon_{i1},\ldots,\varepsilon_{in})'$, $i=1,\ldots,N$, are independent vectors of random errors, $\psi_1,\ldots,\psi_N$ are vectors of individual parameters in $\mathbb{R}^d$, $\beta \in \mathbb{R}^{p}$ is a fixed-effect parameter that is common to the $N$ subjects and the regression function $f$ is known and nonlinear in at least one component of parameter $(\psi_i,\beta)$. In applications where the observations are kinetics or growth curves, it is common for the regression function $f$ to be the solution of an ordinary differential equation (ODE) as shown in the two examples below.
\item at second level, the individual parameters $\psi_1,\ldots,\psi_N$ are defined as independent and identically distributed ({\it i.i.d.}) random variables with common density function $g$. In the most common uses of nonlinear mixed-effects models, the distribution of the random effects is specified up to an unknown parameter $\tau$, \textit{i.e.} 
\begin{equation*}
\psi_i \underset{i.i.d.}{\sim} g(\psi,\tau)d\psi \, , \, i=1,\ldots,N.
\end{equation*}
The value of $\tau$ characterizes the typical characteristics of the phenomenon under study and the inter-individual variability in the population of interest.
\end{enumerate}

In many recent applications, the need to account for random fluctuations in the measured dynamics unrelated to measurement errors, but possibly caused by unknown biological processes, has led to replace the nonlinear regression function $f$ by a continuous-time stochastic process. This has given rise to nonlinear mixed-effects models based on stochastic differential equations where
\begin{align}
Y_{ij} & =  X_i(t_{ij}) + \varepsilon_{ij} \, , \, i=1,\ldots,N \, , \, j=1,\ldots,n, \nonumber\\
dX_i(t) & = b(X_i(t),\psi_i,\beta,t) dt + \sigma(X_i(t),\psi_i,\beta,t) dW_i(t), \, X_i(0) = x_{i0}, \nonumber\\
\psi_i & \underset{i.i.d.}{\sim} g(\psi,\tau) d\psi, \label{eq:noisySDE}
\end{align}
where $(W_i(t),t \geq 0)$, $i=1,\ldots,N$, are independent standard Brownian motions, $x_{i0}$, $i=1,\ldots,N$, are initial conditions for the $N$ processes, $\varepsilon_i=(\varepsilon_{i1},\ldots,\varepsilon_{in})'$, $i=1,\ldots,N$, are independent vectors of random errors, $\psi_1,\ldots,\psi_N$ are vectors of random individual parameters and $\beta$ is a vector of fixed-effects.\\

\textbf{Example. Growth curves.} \cite{lindstrom1990}, \cite{picchini2011} and many other references present the study of the growth of orange trees based on measures of the trunk circumferences of $N=5$ trees at $n=7$ occasions. The model on which there is consensus among the various publications is a logistic growth model:
\begin{align}
Y_{ij} & = \frac{\psi_{i1}}{1+\exp(-(t_{ij}-\psi_{i2})/\psi_{i3})} + \varepsilon_{ij} \, , \, \forall i,j \label{eq:orange.NLME}\\
\psi_i & \underset{i.i.d.}{\sim} \mathcal{N}(\mu,\Omega), \nonumber
\end{align}
where $\psi_i  = (\psi_{i1},\psi_{i2},\psi_{i3})'$. In this example, $\tau=(\mu,\Omega)$ but note that if one component of $\psi_i$ is degenerated, \textit{i.e.} $\psi_{ij} = \mu_j$ $\forall i\in 1,\ldots,N$, then the corresponding $\mu_j$ enters the fixed-effects parameter that was denoted $\beta$ in \eqref{eq:nonlinear}. The nonlinear regression function $$f(\psi,t) = \frac{\psi_1}{1+\exp(-(t-\psi_2)/\psi_3)},$$ is the solution of the ODE $$f'(t) = \frac{1}{\psi_1\psi_3} f(t)(\psi_1-f(t))dt,$$ with initial condition $\psi_1/(1+\exp(\psi_2/\psi_3)$.
In \cite{picchini2011}, a diffusion version of this model has been proposed by replacing function $f$ in \eqref{eq:orange.NLME} by a stochastic process ruled by the following stochastic differential equation
\begin{equation*}
dX(t) = \frac{1}{\psi_1 \psi_3} X(t) (\psi_1 - X(t))dt + \sigma \sqrt{X(t,\psi_1)} dW(t)
\end{equation*}
to account for unexpected growth rate changes in time. The model becomes
\begin{align*}
Y_{ij} & =  X_i(t_{ij}) + \varepsilon_{ij} \, , \, i=1,\ldots,N \, , \, j=1,\ldots,n, \\
dX_i(t) & = \frac{1}{\psi_{i1} \psi_{i3}} X_i(t) (\psi_{i1} - X_i(t))dt + \sigma \sqrt{X(t,\psi_{i1})} dW_i(t),\\
\psi_i & \underset{i.i.d.}{\sim} g(\psi,\tau) d\psi,
\end{align*}
where $(W_i(t),t \geq 0)$, $i=1,\ldots,N$, are independent Brownian motions and $\sigma$ is an unknown fixed-effect in the diffusion coefficient.\\

\textbf{Example. Pharmacokinetics.} Nonlinear regression and SDEs with mixed-effects are also extremely popular in pharmacokinetics (PK) to describe the evolution of drug concentrations over time based on repeated blood samples from several patients. The analysis of the kinetics of the anti-asthmatic drug Theophylline is a well-known PK study. It is based on serum concentrations measured in $N=12$ patients on $n=11$ time occasions after each patient $i=1,\ldots,N$, has received a dose $\mathrm{D_i}$ of the drug. The reader can, for example, consult the books \cite{pinheiro2009} and \cite{lavielle2014} where the statistical analysis of the Theophylline dataset serves as an illustration of mixed-effects model approaches. This study is classically based on the following standard PK model 
\begin{align*}
Y_{ij} & = \frac{\mathrm{D_i} ka_i ke_i}{Cl_i (ka_i - ke_i)}\left(\exp(-ke_i \ t_{ij}) - \exp(-ka_i \ t_{ij})\right) + \varepsilon_{ij} \, , \, \forall i,j\\
\psi_i & \underset{i.i.d.}{\sim} \mathcal{N}(\mu,\Omega) \nonumber,
\end{align*}
where $\psi_i=(ka_i,ke_i,Cl_i)'$ are individual pharmacokinetic parameters (respectively the absorption constant, the elimination constant and the clearance of the drug), $\tau=(\mu,\Omega)$ and the regression function $$f(\psi,t)=\frac{\mathrm{D} \ ka \ ke }{Cl (ka - ke)}(\exp(-ke \ t) - \exp(-ka \ t))$$ verifies
$$
f'(t)=\frac{\mathrm{D} \ ka \ ke}{Cl} \exp(-ka \ t) - ke\ f(t).
$$
In \cite{donnet2008}, Donnet and Samson introduce the following stochastic differential equations instead of $f$ to account for irregularities in the observed kinetics that are not captured by the ODE model: 
\begin{align*}
Y_{ij} & =  X_i(t_{ij}) + \varepsilon_{ij} \, , \, i=1,\ldots,N \, , \, j=1,\ldots,n, \\
dX_i(t) & = \left(\frac{\mathrm{D}_i ka_i ke_i}{Cl_i} \exp(-ka_i t) - ke_i X_i(t) \right) dt + \gamma dW_i(t),\\
\psi_i & \underset{i.i.d.}{\sim} \mathcal{N}(\mu,\Omega),
\end{align*}
where $(W_i(t),t \geq 0)$, $i=1,\ldots,N$, are independent Brownian motions and  $\gamma$ is an additional unknown fixed-effect.

While it is necessary to include a noise term in standard nonlinear models \eqref{eq:nonlinear}, including mixed-effects in SDEs naturally adds noise in the dynamics. It is therefore natural to first look at the theory of estimators in mixed-effects SDEs without measurement errors.

\subsection{Observations and asymptotic framework}
\label{sec:asymptotics}

In the mixed-effects models stated in equations \eqref{eq:nonlinear} and \eqref{eq:noisySDE}, the fixed-effects $\beta$ and the parameters $\tau$ ruling the random effects distribution together help characterize the typical trend and the inter-individual variability in the population of interest. Therefore, providing good estimations is an important issue for the understanding of the phenomenon which is described through the mixed-effects model. It is also essential to know the general properties of parameter estimators for appropriate use in applications depending on the nature of the data available and on the model parameterization. 
While for $k$ \textit{i.i.d.} observations, the unique natural asymptotic framework in which to study the parameter estimators is that $k$ goes to infinity, observing one or several dynamics over time does not lead to a unique asymptotic situation. More broadly, a quick look at the state of the art on theoretical inference for stochastic processes and nonlinear mixed-effects models shows that the properties of the estimators are expected to differ according to the scheme of observations.  

Observing a single dynamics over time falls within the literature of classical stochastic differential equations. The properties of estimators then depend on the nature of the observations, continuous or discrete, and if discrete, on the time interval of observation $[0,T]$ ($T \rightarrow \infty$ or $T$ fixed) and on the sampling time interval $\Delta$ between consecutive observations ($\Delta \rightarrow 0$ or $\Delta$ fixed) while the number of observations $n \rightarrow \infty$. We refer the reader to the Appendix for a summary of the existing results in these different situations.

Describing several dynamics simultaneously relies on mixed-effects models. In classical nonlinear mixed-effects models, it is implicitly considered that observations are collected in discrete times over a fixed time interval. The different asymptotic frameworks thus involve the number of subjects $N$ and the number of observations per subject $n$. To the best of our knowledge, there are only results concerning the maximum likelihood estimator (MLE) due to Nie in a series of three publications (\cite{nie2005}, \cite{nie2006}, \cite{nie2007}) that cover the three following situations in standard nonlinear mixed-effects models:
\begin{enumerate}
\item \textbf{Case (1):} the number of subjects $N \rightarrow \infty$ while the number of measurements per subject $n$ is finite. In this situation, the MLE for both the fixed-effects $\beta$ and the random-effects distribution parameters $\tau$ are consistent (\cite{nie2006}) and asymptotically gaussian with convergence rate $\sqrt{N}$ (\cite{nie2007}).  
\item \textbf{Case (2):} the number of subjects $N$ is limited ($N < \infty$) - for ethical reasons for instance - but the number of observations per subject may be large ($n \rightarrow \infty$). This is the most delicate situation because it is often the simplest experimental design to consider in practice, but it is also the one for which estimators can have the poorest behavior. Indeed, when $n \rightarrow \infty$ while $N < \infty$, only the fixed-effects $\beta$ can be consistently estimated but the correct estimation of $\tau$ is not guaranteed in theory (\cite{nie2005}, \cite{nie2007}). 
\item \textbf{Case (3):} both the number of subjects and the number of observations per subject tend to infinity: $N \rightarrow \infty$ and $n \rightarrow \infty$. The properties differ from parameter to parameter. The MLE for the fixed-effects $\beta$ is $\sqrt{Nn}$-consistent whereas the MLE for the random-effects parameters $\tau$ is $\sqrt{N}$-consistent (\cite{nie2007}). 
\end{enumerate}

\subsection{Linear stochastic differential equations with mixed-effects}

SDEs with mixed-effects however do not fall within the scope of Nie's works \cite{nie2005}, \cite{nie2006}, \cite{nie2007}. Prior the work of Ditlevsen et al. (2005) \cite{ditlevsen2005}, no studies had been conducted on estimators for SDEs with mixed-effects. Since then, maximum likelihood estimation and minimum contrast estimation have been theoretically investigated within simplified versions of the SDEs with mixed-effects stated in \eqref{eq:noisySDE}. In practice, it often occurs that the individual dynamics are observed with measurement errors (see models \eqref{eq:nonlinear} and \eqref{eq:noisySDE}). To our knowledge, there are no theoretical result available. Therefore, we have rather not addressed this additional problem in this review paper. Combining stochastic differential equations with random effects raises many technical difficulties for inference. Hence attention has been exclusively limited to linear SDEs with mixed-effects with distinct parameters in the drift and in the diffusion coefficient for which explicit statistics can be derived for estimation and that provide the ability to study the asymptotic properties of estimators:
\begin{equation}
(\psi,\beta) = (\varphi,\gamma), \,  b(x,\varphi)=\varphi'b(x), \, b(x) = (b_1(x),\ldots,b_d(x))', \, \sigma(x,\gamma)=\gamma^{-1/2}\sigma(x).
\label{eq:linear}
\end{equation}
 For sake of clarity, we restrict our attention to this class of model and to the case where the SDEs are unidimensional:
\begin{equation}
\label{eq:model}
dX_i(t) =  \phi_i' b(X_i(t))dt + (\Gamma_i)^{-1/2}\sigma(X_i(t))dW_i(t),  \, X_i(0) = x_{i0} , \; i=1,\ldots, N,
\end{equation}
where $(W_i(t),t \geq 0)$, $i=1,\ldots,N$, are independent Brownian motions, $x_{i0}$, $i=1,\ldots,N$, are the initial conditions of the $N$ processes. The components of $(\phi_i,\Gamma_i) \in \mathbb{R}^d \times \mathbb{R}^{+*}$ may be fixed or random, and the random coordinates are assumed to be \textit{i.i.d.} random variables with a common density $g(\cdot,\tau)$. In what follows, the whole unknown parameter (including $\tau$ and the possible fixed-effects) will be denoted by $\theta$ and its true value will be denoted by $\theta_0$. We also make the necessary assumptions to ensure the existence and uniqueness of a strong solution to \eqref{eq:model}.

Depending on whether random effects enter the drift or the diffusion coefficient or both, the strategies to be adopted for inference differ, leading to different asymptotic behaviors of the estimators. The asymptotic properties of the corresponding estimates rely on different schemes of observations according to the approach. In any case, assuming that the number of individual dynamics $N \rightarrow \infty$ is a necessary condition to guarantee the good properties of the estimators, but the frequency of observations per trajectory proves to play an important role for the inference.

\section{Maximum likelihood estimation}

Computing the maximum likelihood estimator in mixed-effects models requires some integration over the random-effects distribution. This is usually not possible in a closed form making intractable the study of the theoretical properties of the maximum likelihood estimator.

\subsection{Continuous time observations of the $N$ paths}

In order to illustrate the problems involved in maximum likelihood approaches, let us detail the case where the processes $(X_i)$ are continuously observed on $[0,T]$, where the diffusion coefficient is known and equal for all paths to $\sigma(x)$: 
\begin{equation*}
dX_i(t) = b(X_i(t), \phi_i) + \sigma(X_i(t))dW_i(t), \; i=1,\ldots, N,
\end{equation*}
and where the random effect $\phi_i$ has a distribution $g(\varphi,\theta)d\varphi$. Then, conditionally on $\phi_i=\varphi$, the likelihood for path $i=1,\ldots,N$, is derived through the Girsanov formula:
\begin{equation*}
\exp \left(\int_0^{T} \frac{ b(X_i(s), \varphi)}{\sigma^2(X_i(s))} dX_i(s)  - \frac{1}{2}\int_0^{T}  \frac{ b^2(X_i(s), \varphi)}{\sigma^2(X_i(s))} ds \right).
\end{equation*}
Therefore, by independence of the $N$ processes, the likelihood of the $N$ paths is given by

\begin{equation}
L_{N,T}(\theta) = \prod_{i=1}^{N}  \int \exp \left(\int_0^{T} \frac{ b(X_i(s), \varphi)}{\sigma^2(X_i(s))} dX_i(s)  - \frac{1}{2}\int_0^{T}  \frac{ b^2(X_i(s), \varphi)}{\sigma^2(X_i(s))} ds \right)  g(\varphi; \theta) d\varphi.
\label{eq:girsanov}
\end{equation}

The integrals in \eqref{eq:girsanov} generally cannot be evaluated explicitly. To derive $L_{N,T}(\theta)$ in an explicit form, and therefore obtain some theoretical results on the estimators, the study should be restricted to linear stochastic differential equations, where the drift function verifies \eqref{eq:linear}, with Gaussian random effects in the drift, \textit{i.e.}
\begin{eqnarray*}
dX_i(t) & = & \phi_i' b(X_i(t)) + \sigma(X_i(t))dW_i(t), \; i=1,\ldots, N,\\
\phi_i & \underset{i.i.d.}{\sim} & g(\varphi,\theta)d\varphi = \mathcal{N}(\mu,\Omega).
\end{eqnarray*}

In this specific case, $\theta=(\mu,\Omega)$ and the expression of $L_{N,T}(\theta)$ involves the stochastic integrals 
\begin{equation}
\label{eq:UiVi}
U_i := U_i(T) = \int_0^{T} \frac{b(X_i(s))}{\sigma^{2}(X_i(s))} dX_i(s), \; V_i := V_i(T) = \int_0^{T} \frac{b(X_i(s)) b'(X_i(s))}{\sigma^{2}(X_i(s))} ds,
\end{equation}
according to the formula
\begin{align}
\label{eq:vraisexact}
L_{N,T}(\theta)  = \prod_{i=1}^{N} & \left(\frac{1}{\sqrt{\operatorname{det}(I_d+V_i \Omega})} \exp\left[-\frac{1}{2}\left(\left(\mu-V_i^{-1}U_i\right)'R_i^{-1}\left(\mu-V_i^{-1}U_i\right)\right.\right.\right.\nonumber\\
&\left. \left. \left. - U_i'V_i^{-1}U_i\right)\right]\right),
\end{align} 
where $R_i^{-1} = (I_d+V_i\Omega)^{-1}$ and $I_d$ is the $d \times d$ identity matrix. In \cite{delattre2013}, Delattre et al studied the continuous-time MLE when the $N$ processes are observed on the same time interval $[0,T]$ with $T$ fixed ($T<\infty$) in a purely {\it i.i.d.} setting, \textit{i.e.} $x_{i0} = x_0$ known, and proved the following theorem. 
\begin{prop} \label{theo1} Let ${\hat \theta}_{N,T}$ be a maximum likelihood estimator defined as any solution of $ L_{N,T} ({\hat \theta}_{N,T})= \sup_{\theta \in \Theta} L_{N,T} (\theta)$. Under some appropriate assumptions stated in \cite{delattre2013}, 
\begin{enumerate}
\item ${\hat \theta}_{N,T}$ converges in probability to $\theta_0$,
\item the maximum likelihood estimator satisfies, as $N$ tends to infinity,
$$\sqrt{N}({\hat \theta}_{N,T} -\theta_0)\rightarrow_{{\mathcal D}} {\mathcal N}(0, {\mathcal I}^{-1}(\theta_0)),
$$
where ${\mathcal I}$ refers to the Fisher information matrix whose exact formula is given in \cite{delattre2013}. 
\end{enumerate}
\end{prop}

Some extensions of Theorem \ref{theo1} to independent non identically distributed observation paths were proposed later on. Ruse \& al. (2020) proved consistency and asymptotic normality of the MLE in the case of multidimensional SDEs with a linear drift including covariate-dependent coefficients \cite{ruse2020}. Maitra \& Bhattacharya (2016) also proved consistency and asymptotic normality of the continuous-time MLE in a non {\it i.i.d.} setting without covariates but with different observation time intervals $[0,T_i]$ and different initial conditions $x_{i0}$ for each path \cite{maitra2016}. Delattre et al. (2016) also proved the consistency of the continuous-time MLE in a less conventional framework of curve classification where the random effects $\phi_i$ in the drift are distributed according to a mixture of Gaussian distributions. 

\subsection{Discrete time observations of the $N$ paths}
\label{sec:discretelik}

Although the model is formulated in continuous time, the data are mainly available at discrete time points in practice.   For sake of clarity, we restrict our attention to the case where the observations for the $N$ dynamics are regular and synchronous on a common time interval $[0,T]$. Each sample path $(X_i(t))$ is therefore observed at $n+1$ discrete time points $(t_j, j=0,\ldots,n)$ with regular sampling interval $\Delta=\Delta_n=T/n$. The sequence of observations for process $i$ is denoted by $x_i=(x_{i0}, x_{i1},\ldots,x_{in})$, where $x_{ij}=X_i(t_{j})$, $i=1,\ldots,N$, $j=1,\ldots,n$.

It is well-known that for discretely observed SDEs the exact likelihood is intractable except in very special cases, since it relies on the transition densities of $(X_i)$ which generally have no closed form. To be general enough, consider the most general model
\begin{equation*}
dX_i(t) = b(\psi_i,X_i(t))dt + \sigma(\psi_i,X_i(t))dW_i(t), 
\end{equation*}
where $\psi_i \sim g(\psi,\theta)d\psi$. Proceeding as above, we get the following generic expression for the likelihood  of the $N$ discretely observed paths
\begin{equation*}
L_{N,n}(\theta) = \prod\limits_{i=1}^{N} L(x_i;\theta) = \prod\limits_{i=1}^{N} \int{p(x_i|\psi)g(\psi,\theta)d\psi},
\label{eq:SDEMElike}
\end{equation*} 
where, by the Markov property of $(X_i(t),t \geq 0)$, the conditional individual likelihoods are given by the product of $n$ transition densities   
\begin{equation*}
p(x_i|\psi) := p(x_i|\psi_i=\psi) = p(x_{i0};\psi) \prod\limits_{j=0}^{n-1} p(x_{i,j+1},\Delta|x_{i,j};\psi).\\
\end{equation*} 

\textbf{Example. Brownian motion with drift.} In the pioneering work of Ditlevsen et al. (2005) \cite{ditlevsen2005}, where the case of a Brownian motion with drift and Gaussian random effects in the drift is considered:
\begin{eqnarray*}
dX_i(t) &=&  (\phi_i - \frac{1}{2}\sigma^2)dt + \sigma dW_i(t) \, , \,  i=1,\ldots,N,\\
\phi_i & \underset{i.i.d.}{\sim} & \mathcal{N}(\mu,\omega^2),
\end{eqnarray*}
the transition densities of the individual paths are Gaussian densities, and the likelihood is explicit:
\begin{eqnarray*}
L_{N,n}(\theta) &= &\frac{(\eta^2)^{N/2}}{(2 \pi \sigma^2)^{Mn/2}} \prod_{i=1}^{N} \frac{1}{\Delta^{n/2} \sqrt{T+\eta^2}} \times \\
&& \exp\left( - \frac{\sum\limits_{i,j} (x_{ij} - x_{i,j-1} - \alpha \Delta)^2/\Delta - \sum\limits_i (x_{in}-x_{i0}-\alpha T)^2 (T+\eta^2)^{-1} }{2\sigma^2}\right),
\end{eqnarray*}
where $\alpha=\mu-\sigma^2/2$ and $\eta^2=\sigma^2/\omega^2$. Note that the MLEs for parameters $\mu$, $\omega^2$ and $\sigma$ are also available in a closed form and their asymptotic variances can be explicitly derived. Ditlevsen et al. (2005) did not study consistency or asymptotic normality but they discussed the asymptotic behavior of the estimators heuristically when $N \rightarrow \infty$. \\

In more general models than the one considered in the previous example, exact maximum likelihood estimation is rarely feasible. For $L_{N,n}(\theta)$ to be known in a closed form, the transition probabilities of the processes must be explicitly known and integration with respect to the random-effects density must be explicitly feasible. This is rarely the case. We now introduce alternative estimation strategies.

\section{Estimators derived from the continuous-time likelihood}
\label{sec:continuous}

The first method considered in the literature consists in doing as if the observations were collected in continuous-time in order to bypass the difficulty of explicitly knowing the transition densities of the processes. The parameter estimates are therefore computed by maximizing the discretized version of the continuous-time observation likelihood \eqref{eq:girsanov}. This idea originates from Delattre \& al. (2013) who dealt with this approach in the case where the diffusion coefficient is known and equal for all paths to $\sigma(x)$ with linear Gaussian random effects in the drift:
\begin{eqnarray}
dX_i(t) & = & \phi_i' b(X_i(t)) dt + \sigma(X_i(t)) dW_i(t) \, , \, i=1,\ldots,N, \label{eq:lineardrift}\\
\phi_i & \underset{i.i.d.}{\sim} & g(\varphi,\theta)d\varphi = \mathcal{N}(\mu,\Omega), \nonumber
\end{eqnarray}
in a purely \textit{i.i.d.} setting where the known initial condition $x_{i0} = x_0$ is the same for all paths and the $n$ observations per path are synchronous on the same fixed time interval $[0,T]$, $T<\infty$. The unknown parameters to be estimated here are $\theta=(\mu,\Omega)$, and the estimators are studied under the asymptotic framework where $N,n \rightarrow \infty$.\\

By replacing the terms $U_i$, $V_i$ and $R_i^{-1}$, $i=1, \ldots, N$, by their discretized versions 
\begin{align}
U_{i,n} & = \sum_{j=0}^{n-1} \frac{b(x_{i,j})}{\sigma^2(x_{i,j})} (x_{i,j+1} -  x_{i,j}),\; V_{i,n} =\sum_{j=0}^{n-1} \frac{b(x_{i,j})b'(x_{i,j})}{\sigma^2(x_{i,j})}(t_{j+1} - t_j), \label{eq:UinVin} \\
R_{i,n}^{-1} & = (I_d+V_{i,n}\Omega)^{-1} \nonumber, 
\end{align}
in the continuous-time likelihood, whose expression is \eqref{eq:vraisexact} in models of the form \eqref{eq:lineardrift}, and after a logarithmic transformation, Delattre et al. (2013) propose the following contrast for inference
\begin{align}
\ell^n_N(\theta) = & -\frac{1}{2}\sum_{i=1}^{N} \log \det (I_d + V_{i,n} \Omega) -\frac{1}{2} \sum_{i=1}^{N} (\mu-V_{i,n}^{-1}U_{i,n})'R_{i,n}^{-1} (\mu-V_{i,n}^{-1}U_{i,n}). \label{eq:continuous.contrast}
\end{align} 
The parameter estimators are thus obtained as
\begin{equation*}
\widehat{\theta}_N^{(n)} = \underset{\theta \in \Theta}{\operatorname{argmin}} \; \ell^n_N(\theta).
\end{equation*}

Let us give some concrete examples. \\

\textbf{Examples.}
\begin{enumerate}
\item \textit{Brownian motion with drift.} $dX_i(t)=\phi_i dt + dW_i(t)$, $\phi_i \sim \mathcal{N}(\mu,\omega^2)$. \\
The exhaustive statistics $U_{i,n}=(x_{i,n}-x_{i0})$ and $V_{i,n}=T$ are explicit. The contrast expression is very simple
$$
\ell^n_N(\theta) = -\frac{N}{2} \log(1+\omega^2T) - \frac{1}{2(1+\omega^2 T)} \sum_{i=1}^{N}\left(\mu-\frac{1}{T}(x_{i,n}-x_{i0})\right)^2, 
$$
and leads to explicit estimators for both parameters: 
$$
\hat{\mu}_N^{(n)} = \frac{1}{NT} \sum_{i=1}^{N} (x_{i,n}-x_{i0}) \, , \, \widehat{\omega^2}_N^{(n)} = \frac{1}{T}\left( \frac{1}{N} \sum_{i=1}^{N} \left(\hat{\mu}_N^{(n)}  - \frac{x_{i,n}-x_{i0}}{T}\right)^2 -1 \right).
$$
Note that in this very simple model, $\ell^n_N(\theta)$ coincides with the continuous-time log-likelihood. Thus the continuous-time MLEs and the discrete estimators for parameters $\mu$ and $\omega^2$ coincide.
\item \textit{Ornstein-Uhlenbeck.} $dX_i(t)=\phi_i X_i(t) dt + dW_i(t)$, $\phi_i \sim \mathcal{N}(\mu,\omega^2)$. \\
The continuous-time pivotal integrals $U_i$ and $V_i$ and their discretized versions are given by $$U_i=\int_0^T X_i(s) dX_i(s), \,V_i=\int_0^T X_i(s) ds,$$ and $$U_{i,n} = \sum\limits_{j=0}^{n-1} x_{i,j}(x_{i,j+1}-x_{i,j}), \, V_{i,n}=\Delta \sum\limits_{j=0}^{n-1} x_{i,j}^2$$ respectively. Here, $\hat{\mu}_N^{(n)}$ and $\widehat{\omega^2}_N^{(n)}$ do not have analytic expressions, but the contrast $\ell^n_N(\theta)$ leads to explicit estimating equations that are easy to solve numerically. 
\end{enumerate}
As shown in both previous examples, the contrast or approximate likelihood defined in \eqref{eq:continuous.contrast} is explicit, which is very interesting because this is rarely the case in mixed-effects models. Moreover, it leads to estimators that are both easy to compute and provide the ability to study their asymptotic properties when the number of paths $N$ and the number of observations per trajectory $n$ tend to infinity. Delattre et al (2013) obtain the following result.
\begin{prop} \label{theo2} Under some appropriate assumptions stated in \cite{delattre2013}, and considering $N\rightarrow \infty$,
\begin{enumerate}
\item if $n\rightarrow +\infty$, then $\widehat{\theta}_N - \widehat{\theta}_N^{(n)} = o(1)$,
\item if $n=n(N)\rightarrow +\infty$ in such a way that $\frac nN \rightarrow +\infty$, then $\sqrt{N} (\widehat{\theta}_N - \widehat{\theta}_N^{(n)})  =o(1)$.
\end{enumerate}
\end{prop}

Good performance of $\widehat{\theta}_N^{(n)}$ is therefore achieved provided that $N,n\rightarrow +\infty$ and $\frac nN \rightarrow +\infty$, and thus, as $T$ is fixed, under the uncommon high-frequency scheme where $\Delta N \rightarrow 0$ while $N\rightarrow\infty$. In other words, the success of the above approach does not only depend on the number of observed dynamics $N$ but also on the sampling interval $\Delta$. If $\Delta$ is not small enough, the approximation of $U_i$ (resp. $V_i$) by $U_{i,n}$ (resp. $V_{i,n}$) may be poor and lead to bias in estimation.
It is interesting to note that contrary to the more usual framework of a purely fixed-effects diffusion discussed in the Appendix section, where one would observe only one trajectory, it is not necessary to place oneself in a long time setting, {\it i.e.} $T \rightarrow +\infty$, to enable consistent estimation of the drift parameters. It is therefore not required to assume here that the diffusions are ergodic or in stationary regime. A second interesting point is that contrast $\ell^n_N(\theta)$ is still valid when $\Omega$ is singular, \textit{i.e.} when there are some fixed (and potentially unknown) components in $\phi_i$. The asymptotic results obtained in \cite{delattre2013} implicitly mean that the fixed effects parameters and the random effects parameters in the drift are estimated with the same rate of convergence. Such a result is unusual in classical mixed-effects models where a higher rate of convergence ($\sqrt{Nn}$) is expected for the fixed-effects parameters when both $N,n \rightarrow \infty$ (see Section \ref{sec:asymptotics} for more details). 
Let us end by stressing that a major limitation of this approach based on the continuous-time likelihood is that it is impossible to build consistent estimators for the diffusion parameters, if any, when $T$ is fixed. 

\section{Contrasts derived from the Euler scheme}

\subsection{General contrast definition}

To deal with more general situations where both the drift parameters and the diffusion parameter are unknown, either fixed or random, a classical way to derive contrasts is to consider the Euler scheme associated with the differential equations. The theoretical aspects of such approach have been investigated in \cite{delattre2015}, \cite{delattre2017} and \cite{delattre2018} in the particular case of linear SDEs, where the drift and the diffusion coefficient verify \eqref{eq:linear}:
\begin{equation*}
dX_i(t) = \phi_i'b(X_i(t))dt + \Gamma_i^{-1/2} \sigma(X_i(t))dW_i(t), \, X_i(0)=x_{i0}, \, i=1,\ldots,N,
\end{equation*}
with appropriate distributions for the random components in $\psi_i=(\phi_i,\Gamma_i)$.

As discussed in Section \ref{sec:discretelik}, the discrete-time likelihood is not available in a closed form. To obtain an explicit contrast, the transition densities $p(x_{i,j+1}, \Delta |x_{i,j};\psi)$ are approximated with the following Gaussian densities:
\begin{equation*}
\tilde{p}(x_{i,j+1}, \Delta| x_{i,j};\psi=(\varphi,\gamma)) = \frac{1}{\sqrt{2\pi \Delta \gamma^{-1}\sigma^2(x_{i,j})}} \exp\left\{-\frac{x_{i,j+1}-x_{i,j}-\Delta \varphi' b(x_{ij})}{2 \Delta \gamma^{-1} \sigma^2(x_{i,j})}\right\}.
\end{equation*}

Then, conditionally on $\psi_i=(\varphi,\gamma)$, the likelihood for path $i=1,\ldots,N$, is approximated by

\begin{eqnarray*}\label{likeuler}
\tilde{p}(x_i|\psi=(\varphi,\gamma)) & = & \prod_{j=0}^{n-1} \tilde{p}(x_{i,j+1}, \Delta |x_{i,j};\psi) \nonumber\\
& = & \gamma^{n/2} \exp{[-\frac{\gamma}{2}(S_{i,n} +\varphi' V_{i,n} \varphi  -2\varphi' U_{i,n})]},
\end{eqnarray*}
where

\begin{equation}
S_{i,n}=  \frac{1}{\Delta}\sum_{j=1}^{n}\frac{\left(x_{i,j}-x_{i,j-1}\right)^2 }
{\sigma^2(x_{i,j-1})}, 
\label{eq:Si}
\end{equation}

and $U_{i,n}$ and $V_{i,n}$ are the same statistics as the ones defined in \eqref{eq:UinVin}. Contrasts for inference are then built by integrating out the conditional likelihood of the Euler schemes with respect to the densities of the random effects. This does only lead to a closed form expression if specific  distributions are chosen for the random effects. In \cite{delattre2015} and  \cite{delattre2017}, the author consider Gamma distributed random effects in the diffusion coefficient or Gaussian random effects in the drift, leading in both cases to an explicit expression for the integrated densities of the Euler schemes associated with the $N$ paths. Combining random effects in the drift and random effects in the diffusion coefficient is nevertheless not straightforward. In \cite{delattre2018}, bringing together Gamma random effects in the diffusion coefficient and conditionally Gaussian random effects in the drift, and carrying out appropriate modifications to the integrated densities of the Euler schemes, the authors provide two explicit contrasts. Except for one very specific situation where the contrasts derived from the Euler schemes associated to the paths coincide with the exact likelihood of the processes, high-frequency observations are required for them to provide consistent estimators. The respective proofs in  \cite{delattre2015}, \cite{delattre2017} and \cite{delattre2018} are developed in a purely \textit{i.i.d.} setting where $x_{i0}= x$, $i=1,\ldots,N$ and observations are collected in a finite time interval \textit{i.e.} $T<\infty$. We now detail the results for the three different cases.\\

\subsection{Fixed effect in the drift and random effect in the diffusion coefficient.} 
\label{sec:euler.rediff}

The situation where the model does only include random effects in the diffusion coefficient has been addressed in \cite{delattre2015} and \cite{delattre2017}. In \cite{delattre2015}, the drift is set to 0 or ignored:
\begin{equation*}
dX_i(t) = \Gamma_i^{-1/2} \sigma(X_i(t)) dW_i(t),
\label{eq:nuldrift}
\end{equation*}
whereas \cite{delattre2017} deals with the more general case of a non-null drift that includes a linear unknown fixed parameter:
\begin{equation}
dX_i(t) = \varphi'b(X_i(t))dt + \Gamma_i^{-1/2} \sigma(X_i(t)) dW_i(t).
\label{eq:nonnuldrift}
\end{equation}
In both cases, the author consider Gamma distributed random effects in the drift, 
\begin{equation}
\Gamma_i \underset{i.i.d.}{\sim} g(\gamma,\tau)d\gamma = G(a,\lambda),
\label{eq:gamma.re}
\end{equation} 
thus the parameters to be estimated are $\theta=(\lambda, a)$ if the drift is set to 0 and $\theta=(\lambda, a, \varphi)$ otherwise. \\

When the drift is non-null (model defined by equations \eqref{eq:nonnuldrift}-\eqref{eq:gamma.re}), the pseudo-log-likelihood derived from the Euler schemes of the $N$ paths is: 

\begin{eqnarray}
\tilde{\ell}_{N,n}(\theta) & = & N a \log(\lambda) + N \log(\Gamma(a+n/2)) - N \log(\Gamma(a)) \nonumber\\
& & - (a+n/2) \sum_{i=1}^{N} \log \left(\lambda+\frac{1}{2} (S_{i,n} - 2 \varphi U_{i,n} + \varphi^2 V_{i,n} )\right),
\label{eq:euler.gamma}
\end{eqnarray}
where $U_{i,n}$, $V_{i,n}$ and $S_{i,n}$ are respectively defined in \eqref{eq:UinVin} and \eqref{eq:Si}, and $\Gamma(z)$ is the Gamma function. Unsurprisingly, setting $\varphi=0$ in \eqref{eq:euler.gamma} displays the pseudo-log-likelihood studied in \cite{delattre2015} for the particular model without drift.

The following estimator ${\tilde \theta}_{N,n}$ is then proposed

\begin{equation}
\nabla_{\theta} \tilde{\ell}_{N,n}({\tilde \theta}_{N,n}) = 0,
\label{eq:estim.euler.re.diff}
\end{equation}
where $\nabla_{\theta} \tilde{\ell}_{N,n}(\theta)$ is the pseudo-score associated to the contrast defined in \eqref{eq:euler.gamma}.

Contrast $\tilde{\ell}_{N,n}(\theta)$ leads to estimators that are weakly consistent and asymptotically Gaussian with unique rate of convergence $\sqrt{N}$ for all components of ${\widetilde \theta}_{N,n}$ when both $N,n \rightarrow \infty$. 

\begin{prop} \label{theo3}
Under appropriate assumptions stated in \cite{delattre2017}, if $N,n$ tend to infinity with $N/n \rightarrow 0$ (\textit{i.e.} $N \Delta \rightarrow 0$), then 
\begin{enumerate}
\item a solution ${\widetilde \theta}_{N,n}=({\widetilde \lambda}_{N,n}, {\widetilde a}_{N,n}, {\widetilde \varphi}_{N,n})$ to \eqref{eq:estim.euler.re.diff} exists with probability tending to $1$ which is consistent, 
\item this solution ${\widetilde \theta}_{N,n}$ is such that $\sqrt{N}({\widetilde \theta}_{N,n}-\theta_0)$ converges in distribution to ${\mathcal N}(0, {\mathcal V}^{-1}(\theta_0))$  where 
\begin{equation*}
{\mathcal V}(\theta)= \left(\begin{array}{c|c}
I_0(\lambda,a) & \textbf{0}\\
\hline 
\textbf{0} & {\mathbb E}_{\theta } \left(\Gamma V(T) \right)
\end{array}\right),
\end{equation*}
$I_0(\lambda,a)$ is the Fisher information matrix of an \textit{i.i.d.} sample of a Gamma distribution $G(a,\lambda)$, $\Gamma$ is a $G(a,\lambda)$ distributed random variable and $V(T)$ is defined in \eqref{eq:UiVi}.
\end{enumerate}
For the components ${\widetilde a}_{N,n}$ and ${\widetilde \lambda}_{N,n}$ of ${\widetilde \theta}_{N,n}$, the constraints on $N$ and $n$ can be relaxed and the same results can be obtained assuming only $N/n^2\rightarrow 0$ (\textit{i.e. $N \Delta^2 \rightarrow 0$}). 
\end{prop}

It is interesting to note that for parameters $a$ and $\lambda$, provided that $N/n^2 \rightarrow +\infty$, high frequency sampling of the trajectories $(X_i)$ yields an estimator that is asymptotically equivalent to the exact MLE based on direct observation of the $\Gamma_i$. Contrary to the case of a fixed diffusion parameter that is detailed right after, all parameters are estimated here with the same rate $\sqrt{N}$. In \cite{delattre2015} where the authors considered the special case where $\varphi=0$ with the same contrast, the same results were obtained under the stronger constraint $N/n\rightarrow 0$. 

Finally let us note that if $\sigma(\cdot) \equiv 1$, the likelihood of the Euler scheme of the individual processes is the exact likelihood. In this case, $n$ may be finite meaning that high-frequency observations are not required for the exact MLE of $(a,\lambda)$ to be consistent and asymptotically Gaussian (see \cite{delattre2015} for more details). 

\subsection{Random effect in the drift and fixed effect in the diffusion coefficient.} 
\label{sec:euler.redrift}

The case where linear Gaussian random effects in the drift are combined with a linear fixed-effect in the diffusion coefficient
\begin{eqnarray*}
dX_i(t) &=&  \phi_i' b(X_i(t)) dt + \gamma^{-1/2}\sigma (X_i(t)) dW_i(t) \, , \,  i=1,\ldots,N,\\
\phi_i & \underset{i.i.d.}{\sim} & g(\varphi,\tau) d\varphi = \mathcal{N}(\mu,\Omega),
\end{eqnarray*}

has been adressed by Delattre et al (2017) in \cite{delattre2017}. Here, $\theta=(\gamma,\mu,\Omega)$. Integration of the Euler scheme conditional likelihood leads to the following pseudo-log-likelihood: 
\begin{align}
\tilde{\ell}_{N,n}(\theta) = &   \frac{Nn}{2} \log{\gamma}- \frac{1}{2}\sum _{i=1}^N \log (I_d+V_{i,n}\Omega)
 - \frac{\gamma}{2}\sum _{i=1}^N \left[S_{i,n}-U_{i,n}^{-1}V_{i,n}U_{i,n}\right.\\
& \left.+\left(\mu - V_{i,n}^{-1}U_{i,n}\right)'\left((I_d+V_{i,n}\Omega)^{-1} V_{i,n}\right)\left(\mu - V_{i,n}^{-1}U_{i,n}\right)\right],
\label{eq:lNnmult}
\end{align}

and to the following estimator
\begin{equation}
\nabla_{\theta} \tilde{\ell}_{N,n}({\tilde \theta}_{N,n}) = 0,
\label{eq:estim.euler.re.drift}
\end{equation}
where $\nabla_{\theta} \tilde{\ell}_{N,n}(\theta)$ is the pseudo-score associated to $\tilde{\ell}_{N,n}(\theta)$.

Provided some regularity conditions on the drift and the diffusion coefficients, the authors prove that ${\tilde \theta}_{N,n}$ is weakly consistent and asymptotically Gaussian with different rates of convergence for the drift parameters and for the diffusion parameters:
\begin{prop} \label{theo4} Under appropriate assumptions stated in \cite{delattre2017}, if $N,n$ tend to infinity with $N/n \rightarrow 0$, with probability tending to $1$, 
\begin{enumerate}
\item there exists a solution to \eqref{eq:estim.euler.re.drift}, ${\tilde \theta}_{N,n}=(\tilde{\gamma}_{N,n},\tilde{\mu}_{N,n}, \tilde{\Omega}_{N,n})$, which is consistent,
\item this solution is such that
$$D_{N,n}^{-1}({\tilde \theta}_{N,n}-\theta_0) \rightarrow {\mathcal N}(0, {\mathcal I}^{-1}(\theta_0)).$$ 
${\mathcal I}(\theta)$ is the Fisher information matrix and the convergence rate matrix is given by $D_{N,n}=\mathrm{diag}(\frac{1}{\sqrt{Nn}},\frac{1}{\sqrt {N}} I_{d},\frac{1}{\sqrt {N}} I_{d\times d})$. The exact formula of ${\mathcal I}(\theta)$ is provided in \cite{delattre2017}.

\end{enumerate}
\end{prop}

Under the condition $N,n \rightarrow \infty$, $N/n \rightarrow 0$ (\textit{i.e. $N \Delta \rightarrow 0$}), the rate of convergence is $\sqrt{N}$ for the random effects parameters, whereas the deterministic drift component $\gamma$ has a faster rate $\sqrt{Nn}$. If there were only fixed effects in the drift in addition to a fixed parameter in the diffusion coefficient, contrast \eqref{eq:lNnmult} would still be valid since $\Omega$ does not need to be invertible. Parameter $\mu$ would then still be estimated with the same rate $\sqrt{N}$. This differs from standard nonlinear mixed effects models where, according to Nie's results in \cite{nie2007} (see Section \ref{sec:asymptotics}), a faster rate would be expected.

\subsection{Random effects in the drift and in the diffusion coefficient.}
\label{sec:euler.redriftdiff}

The latter situation is studied in \cite{delattre2018}. The random effects enter linearly both the drift and the diffusion coefficient. The random effects in the diffusion coefficient follow a Gamma distribution whereas the random effects in the drift are conditionally Gaussian:
\begin{eqnarray*}
dX_i(t) &=&  \phi_i' b(X_i(t)) dt + \Gamma_i^{-1/2}\sigma (X_i(t)) dW_i(t) \, , \,  i=1,\ldots,N,\\
\Gamma_i &\underset{i.i.d.}{\sim} & G(a,\lambda), \; \phi_i| \Gamma_i = \gamma \sim \mathcal{N}(\mu,\gamma^{-1}\Omega).
\end{eqnarray*}
The parameter to be estimated is thus $\theta=(\lambda,a, \mu, \Omega)$. In this specific situation, the integral of the conditional likelihood of the Euler schemes of the processes $(X_i)$ over the random effects distribution is not always defined thus raising technical difficulties.  Two alternative contrasts are therefore proposed by the authors. \\

\textit{Contrast 1.} The first contrast is obtained by truncating the integral of the conditional likelihood of the Euler schemes over the random effects distribution on a subset where it is defined:
\begin{eqnarray}
\textbf{U}_{N,n}(\theta)& = & N \left[a\log \lambda-\log(\Gamma(a)) +\log(\Gamma(a+(n/2)))-(a+(n/2))\log(a+(n/2))\right] \nonumber\\
&& -\frac{1}{2}  \sum_{i=1}^{N}\log{(I_d+ V_{i,n} \Omega)}  - (a+(n/2)) \sum_{i=1}^{N} \mathds{1}_{F_{i,n}}\log (Z_{i,n}), 
\label{eq:contrast.u}
\end{eqnarray}
with
\begin{eqnarray*}
Z_{i,n} & = & Z_{i,n}(\theta)= \frac{ 2 \lambda +  S_{i,n}+ (\mu - V_{i,n}^{-1}U_{i,n})'R_{i,n}^{-1}(\mu - V_{i,n}^{-1}U_{i,n}) - U_{i,n}'V_{i,n}^{-1}U_{i,n}}{2 a+n}, 
\end{eqnarray*}
$S_{i,n}$, $U_{i,n}$, $V_{i,n}$ and $R_{i,n}$ are defined in \eqref{eq:Si} and \eqref{eq:UinVin} respectively, and
\begin{equation*}
F_{i,n}=\{S_{i,n}- M_{i,n} \ge \alpha \sqrt{n}\} \, , \, M_{i,n}  = \max\{ c_1+2,2 m^2 \} ( 1+|U_{i,n}|^2 ),
\end{equation*} 
where $c_1$ and $m$ are upper bounds for $||\mu||$ and the maximum eigenvalue of $\Omega$ respectively, and $\alpha >0$ is a constant to be defined. 

This leads to the definition of a new estimator
\begin{equation}
\nabla_{\theta} \textbf{U}_{N,n}({\tilde \theta}_{N,n}^U) = 0.
\label{eq:estim.contrast1}
\end{equation}

The following result is derived. 

\begin{prop}\label{theo5} Under appropriate assumptions stated in \cite{delattre2018}, if $n$ and $N$  tend to infinity with $N/n \rightarrow 0$ then
\begin{enumerate}
\item with probability tending to $1$, a solution ${\tilde \theta}_{N,n}^U$ to \eqref{eq:estim.contrast1} exists which is consistent,
\item it is such that ${\sqrt{N}}({\tilde \theta}_{N,n}^U - \theta_0)$ converges in distribution to ${\mathcal N}(0, {\mathcal J}^{-1}(\theta_0))$ where 
\begin{equation*}
{\mathcal J}(\theta)= \left(\begin{array}{c|c}
 I_0(\lambda,a) & \textbf{0}\\
 \hline 
\textbf{0} & J(\theta)
\end{array}\right),
 \end{equation*}
 $I_0(\lambda,a)$ is the Fisher information matrix of an \textit{i.i.d.} sample of a Gamma distribution $G(a,\lambda)$. The expression of $J(\theta)$ is given in \cite{delattre2018}.
\end{enumerate}
For the first two components, the constraint $N/n^2 \rightarrow 0$ is enough.
\end{prop}

\textit{Contrast 2.}
The second contrast is the sum of two terms, the first one depending on the Gamma distribution parameters and the second one depending on the Gaussian distribution parameters:

\begin{equation}\label{eq:contrast.v}
{\boldsymbol V}_{N,n}(\theta)= {\boldsymbol V}_{N,n}^{(1)}(\lambda,a)+{\boldsymbol V}_{N,n}^{(2)}(\mu, \Omega),
\end{equation}

where

\begin{eqnarray*}
{\boldsymbol V}_{N,n}^{(1)}(\lambda,a) & = &  - (a+\frac{n}{2}) \sum_{i=1}^{N} \log{\left(\lambda + \frac{S_{i,n}}{2}\right)} \\
&& + N \left\{a\log{\lambda} + \log{ \Gamma(a+(n/2)) - \log{\Gamma(a) } }\right\}, \\
{\boldsymbol V}_{N,n}^{(2)}(\mu, \Omega)& = & - \sum_{i=1}^{N}\left\{\frac{n}{2 \,S_{i,n}} \mathds{1}_{S_{i,n}\geq k\sqrt{n}}\left(\mu-V_{i,n}^{-1}U_{i,n}\right)' R_{i,n}^{-1} \left(\mu-V_{i,n}^{-1}U_{i,n}\right) \right.\\
&& \left.+ \frac{1}{2}\log{(\operatorname{det}(I_d+ V_{i,n}\Omega)})\right\},\\
\end{eqnarray*}

and the constant $k$ in the truncations is tuned by the user. 

This leads to a new estimator
\begin{equation}
\nabla_{\theta} \textbf{V}_{N,n}({\tilde \theta}_{N,n}^V) = 0.
\label{eq:estim.contrast2}
\end{equation}

This decomposition of ${\boldsymbol V}_{N,n}(\theta)$ as the sum of ${\boldsymbol V}_{N,n}^{(1)}$ and ${\boldsymbol V}_{N,n}^{(2)}$ allows the drift parameters $(\mu,\Omega)$ and the diffusion parameters $(\lambda,a)$ to be estimated separately. Compared to ${\boldsymbol U}_{N,n}$, which is difficult to implement due to its dependence on the parameter space $\Theta$ which is unknown in practice, ${\boldsymbol V}_{N,n}$ is easy to use and leads to a numerically more stable estimation procedure. The authors show that the estimators provided by ${\boldsymbol U}_{N,n}$ and ${\boldsymbol V}_{N,n}$ are asymptotically equivalent.

\begin{prop} \label{theo6} Under appropriate assumptions stated in \cite{delattre2018}, if $n$ and $N$  tend to infinity with $N/n \rightarrow 0$ (\textit{i.e.} $N \Delta  \rightarrow 0$) then
\begin{enumerate}
\item with probability tending to $1$, a solution ${\tilde \theta}_{N,n}^V$ to \eqref{eq:estim.contrast2} exists which is consistent;
\item it is such that ${\sqrt{N}}( {\tilde \theta}_{N,n}^V - \theta_0)$ converges in distribution to ${\mathcal N}(0, {\mathcal J}^{-1}(\theta_0))$ where ${\mathcal J}(\theta)$ is defined in Theorem \ref{theo5}.
\end{enumerate}
For the first two components, the constraint $N/n^2 \rightarrow 0$ (\textit{i.e.} $N \Delta^2 \rightarrow 0$) is enough.
\begin{enumerate}
\setcounter{enumi}{2}
\item The estimators ${\tilde \theta}_{N,n}^V$ and ${\tilde \theta}_{N,n}^U$ are asymptotically equivalent.
\end{enumerate}
\end{prop}  

It is worth noting that with both contrasts, the estimators of $a$ and $\lambda$ are asymptotically equivalent to the exact maximum likelihood of the same parameters based on the direct observation of $(\Gamma_i,i=1,\ldots,N)$, meaning that there is no loss of information for the diffusion parameters. For the parameters $(\mu,\Omega)$, the constraint $N/n\rightarrow 0$ cannot be weakened and there is a loss of information with respect to the direct observation of the random effects.
 
\subsection{Examples and remarks}

Whatever the situation among those described in Sections \ref{sec:euler.rediff}, \ref{sec:euler.redrift} and \ref{sec:euler.redriftdiff}, the only central quantities to compute the parameter estimators are the pivotal statistics $U_{i,n}$, $V_{i,n}$ and $S_{i,n}$. These statistics are explicit. Consider for example a simple Brownian motion with drift
$$dX_i(t) = \phi_i dt + \Gamma_i^{-1/2}dW_i(t),$$
whether combined with \textit{i)} $\phi_i = \mu$ and $\Gamma_i \sim G(a,\lambda)$, \textit{ii)} $\phi_i\sim \mathcal{N}(\mu,\omega^2)$ and $\Gamma_i = \gamma$ or \textit{iii)} $\Gamma_i \sim G(a,\lambda)$ and $\phi_i|\Gamma_i=\gamma \sim \mathcal{N}(\mu,\gamma^{-1}\omega^2)$. Then, these statistics are given by $U_{i,n}=(x_{i,n}-x_{i0})$, $V_{i,n}=T$ and $S_{i,n}=\sum\limits_{j=1}^{n}(x_{i,j}-x_{i,j-1})^2$. Although the estimators for $\mu$, $\omega^2$, $a$ and $\lambda$ are not explicit, the estimating equations have an analytical form in any of the three situations. Moreover, according to the Theorems \ref{theo4} to \ref{theo6} stated above, provided that $N,n\rightarrow \infty$ with the appropriate constraints between $N$ and $n$ ($N/n \rightarrow 0$ or $N/n^2 \rightarrow 0$) the estimators are $\sqrt{N}$-consistent whatever the situation, except for parameter $\gamma$ when it is a fixed-effect in the diffusion coefficient (situation \textit{ii)}) which has a faster rate $\sqrt{Nn}$.

Let us finally add that the contrast estimation methods presented in Sections \ref{sec:euler.rediff}, \ref{sec:euler.redrift} and \ref{sec:euler.redriftdiff} have been implemented in a \texttt{R}-package \texttt{MsdeParEst} \cite{msdeparest2017}.\\

In the following two sections, we discuss two-step methods, either parametric (see \cite{delattre2015}, Section \ref{sec:plugin.parametric}) or non parametric (see \cite{dion2016a} and \cite{dion2016b}, Section \ref{sec:nonparametric}), whose common intuitive central idea is starting with computing appropriate estimations $\hat{\phi}_1,\ldots,\hat{\phi}_N$ (resp. $\hat{\Gamma}_1,\ldots,\hat{\Gamma}_N$) of the random effects $\phi_1,\ldots,\phi_N$ (resp. $\Gamma_1,\ldots,\Gamma_N$) based on the discrete observations of the $N$ trajectories, and then proceeding classical inference technique as if these estimations were true observations of the random effects.

\section{Parametric estimation based on estimators of the random effects}
\label{sec:plugin.parametric}

In \cite{delattre2015}, the author consider a linear random effect in the diffusion coefficient (see equation \eqref{eq:linear}) and a drift set to 0, {\it i.e.} model \eqref{eq:model} in the particular case where:
\begin{eqnarray*}
dX_i(t) & = & (\Gamma_i)^{-1/2} \sigma(X_i(t)) dW_i(t), \; X_i(0) = x_0, i=1,\ldots, N,\\
\Gamma_i & \underset{i.i.d.}{\sim} & G(a,\lambda).
\end{eqnarray*}
The initial condition $x_0$ is fixed and known. The parameter to be estimated is $\theta=(a,\lambda)$ and the observations available for estimation are obtained over a finite time interval ($T<\infty$). To estimate $\theta$, they propose to estimate the random effects in the diffusion coefficient with the quadratic variation estimator, and to replace the random variables $\Gamma_i$ by their estimators in the likelihood of $(\Gamma_1,\ldots,\Gamma_N)$:
\begin{equation*}
\ell_{\Gamma,N}(\theta) = Na\log\lambda-N\log\Gamma(a) + (a-1)\sum_{i=1}^{N} \log(\Gamma_i) - \lambda \sum_{i=1}^{N} \Gamma_i,
\end{equation*}
where $\Gamma(z)$ is the Gamma function. The usual quadratic variation estimator is however not completely appropriate in the general case where $\sigma(\cdot) \not\equiv 1$ due to moment properties, thus truncated: 
\begin{equation*}
\tilde{\Gamma}_i = \frac{n}{S_{i,n}} \mathds{1}_{(S_{i,n}/n \geq k/\sqrt{n})}, \; \widetilde{\log}  \Gamma_i = \log \left(\frac{n}{S_{i,n}} \right) \mathds{1}_{(S_{i,n}/n \geq k/\sqrt{n})},
\end{equation*}
with $k$ a constant to be tuned by the user. This leads to the following contrast function

\begin{equation}
\tilde{\ell}_{\Gamma,N,n}(\theta) = N a \log \lambda - N \log \Gamma(a) + (a-1) \sum_{i=1}^{n} \widetilde{\log}  \Gamma_i  - \lambda \sum_{i=1}^{N} \tilde{\Gamma}_i, 
\label{eq:contrast.plugin}
\end{equation}
and to the estimator
\begin{equation}
\hat{\theta}_{\Gamma,N,n} = \underset{\theta \in \Theta}{\operatorname{argmin}} \; \tilde{\ell}_{\Gamma,N,n}(\theta).
\label{eq:plugin}
\end{equation}

In \cite{delattre2015}, the drift is set to 0 to study the respective roles of $N$ and $n$ in the asymptotic properties of $\hat{a}_{\Gamma,N,n}$ and $\hat{\lambda}_{\Gamma,N,n}$. The following result is derived.

\begin{prop} \label{theo7}
Under appropriate assumptions stated in \cite{delattre2015}, if $N,n$ tend to infinity in such a way that $\sqrt{N}/n$ tends, then 
\begin{enumerate}
\item an estimator $\hat{\theta}_{\Gamma,N,n}$ which solves \eqref{eq:plugin} exists with probability tending to $1$ and is weakly consistent, \item moreover,  $\sqrt{N}(\hat{\theta}_{\Gamma,N,n}-\theta_0)$ converges in distribution to ${\mathcal N}(0, {\mathcal I}^{-1}(\theta_0))$ where ${\mathcal I}(\theta)$ is the Fisher information matrix of an \textit{i.i.d.} sample of a Gamma distribution $G(a,\lambda)$.
\end{enumerate}
\end{prop}

In other words, provided that $\sqrt{N}/n \rightarrow 0$, the plug-in estimation method causes no loss of information with respect to the case where the random effects would be directly observed. It should be mentioned that if $\sigma(\cdot) \equiv 1$, the contrast expression \eqref{eq:contrast.plugin} is simpler since the truncation is not required. Finally, it is worth noting that the plug-in estimation strategy and the asymptotic properties of the corresponding estimates are also valid in the case of a nonnul drift. 

\textit{Remark:} The plug-in estimation strategy introduced in \cite{delattre2015} for Gamma distributed random effects in the diffusion coefficient could be easily extended to other distributions than a Gamma distribution provided appropriate truncations are introduced. We could imagine reproducing the approach to estimate the drift population parameters based on estimators of the drift random effects since $U_{i,n} V_{i,n}^{-1}$ is a natural estimator for $\phi_i$. 

Genon-Catalot and Lar\'edo (2016, \cite{gecatalot2016}) investigated an analogous plug-in technique in the Ornstein-Uhlenbeck model with one multiplicative Gamma distributed random effect in the drift:
\begin{eqnarray*}
dX_i(t) & = & -\phi_i X_i(t) dt + dW_i(t), \; X_i(0) = \eta_i, \\
\phi_i & \underset{i.i.d.}{\sim} & G(a,\lambda), \; i=1,\ldots,N.
\end{eqnarray*}
The initial conditions $(\eta_1,\ldots,\eta_N)$ are random, and the processes $(X_i(t))$, $i=1,\ldots,N$ are continuously observed on the time interval $[0,T]$ as $T\rightarrow \infty$. If $N,T \rightarrow \infty$ and $N/T \rightarrow 0$, then $\theta=(a,\lambda)$ is estimated  at the rate $\sqrt{N}$. Discrete observations are however not discussed.

\section{Nonparametric estimation of the density of the random effects}
\label{sec:nonparametric}

Let us finish the review by saying a few words about nonparametric methods that have also been investigated for the estimation of the random effects distribution in mixed-effects diffusions. The contributions are fewer in number than those on parametric estimation and concern less general models, including random effects in the drift only and at most bivariate random effects (\cite{comte2013}, \cite{dion2016a} and \cite{dion2016b}). Some of them are implemented in the \texttt{R} package \texttt{mixedsde} \cite{mixedsde2018}. The nonparametric density estimates are built on continuous-time estimators of the random effects. The latter need to be discretized in practice and the influence of the discretization step on the density estimator properties is naturally addressed in the corresponding papers.

For the sake of simplicity, the underlying ideas of the proposed nonparametric estimators are presented in the simplest cases that were investigated by Comte \& al. (2013). The extensions and improvements that came next are discussed in a second step. Two situations should be distinguished as they do not allow for the same developments:
\begin{enumerate}[label=(\alph*)]
\item Multiplicative random effect in the drift
\begin{eqnarray}
dX_i(t) & = & \phi_i b(X_i(t))dt + \sigma (X_i(t)) dW_i(t), \label{eq:mult.model.nonp}
\end{eqnarray}
\item Linear random effect in the drift
\begin{eqnarray}
dX_i(t) & =  & (\phi_i + b(X_i(t))) dt + \sigma (X_i(t)) dW_i(t), \label{eq:lin.model.nonp} 
\end{eqnarray}
\end{enumerate}
where $(W_1,\ldots,W_N)$ are independent standard Brownian motions, $\phi_1, \ldots,\phi_N$ are {\it i.i.d.} random variables in $\mathbb{R}$ with common unknown density $g(\cdot)$ such that $(\phi_1, \ldots,\phi_N)$ and $(W_1,\ldots,W_N)$ are independent. 

In model \eqref{eq:mult.model.nonp} where the random effect is multiplicative in the drift, it is possible to estimate the random effects trajectory by trajectory by using the continuous-time MLE expression of $\phi_i$ when $\phi_i=\varphi$ is deterministic:
\begin{equation}
\hat{\phi}_{i,T}^{(m)} = V_i(T)^{-1} U_i(T),
\label{eq:phiest1}
\end{equation}
where the statistics $U_i(T)$ and $V_i(T)$ are the same statistics that have already been defined in \eqref{eq:UiVi}. According to some classical results on diffusions, $\hat{\phi}_{i,T}^{(m)}$ is a consistent estimate of $\phi_i$ when the time horizon $T$ tends to infinity. Note that if $b(\cdot) = \sigma(\cdot)$, $\hat{\phi}_{i,T}^{(m)} = \phi_i + T^{-1}(W_i(T)-W_i(0))$, $\hat{\phi}_{i,T}^{(m)}$ estimates $\phi_i$ up to an additive Gaussian noise.

In model \eqref{eq:lin.model.nonp} where the random effect is linear in the drift, the following estimator can be proposed for the random effects
\begin{equation}
\hat{\phi}_{i,T}^{(l)} = T^{-1}\left(X_i(T) - X_i(0) - \int_0^{T}b(X_i(s)ds\right),
\label{eq:phiest2}
\end{equation} 
that is, by simple calculations, $\hat{\phi}_{i,T}^{(l)} = \phi_i + T^{-1}\int_0^T \sigma(X_i(s))dW_i(s)$.
As for $\hat{\phi}_{i,T}^{(m)}$, good properties of $\hat{\phi}_{i,T}^{(l)}$ are ensured when $T \rightarrow +\infty$. 

Kernel methods are generic tools for non-parametric estimation of unknown densities. The idea that is developed in \cite{comte2013}, and then taken up in \cite{dion2016a} and \cite{dion2016b}, is that if one is able to compute estimates $\hat{\phi}_{1,T},\ldots,\hat{\phi}_{N,T}$ of the random effects, it is fairly easy and natural to use kernel methods on these estimates $\hat{\phi}_{1,T},\ldots,\hat{\phi}_{N,T}$ to obtain some nonparametric estimation $\hat{g}_h(\cdot)$ of the density function $g(\cdot)$:
\begin{equation*}
\hat{g}_h(x) = \frac{1}{N} \sum_{i=1}^{N} K_h\left(x-\hat{\phi}_{i,T}\right), 
\end{equation*}
where $ K_h(x) = h^{-1} K(h^{-1}x)$, $K$ is a kernel function and $h>0$ is a smoothing parameter, also called bandwidth. The control of the $L_2$-risk is then unusual due to the fact that estimates are used as input data. Comte \& al. (2013) first investigated it based on continuous time observations of the processes. Their results are interesting because they reveal that contrary to most of the available parametric methods discussed in the previous section, kernel methods require the processes to be observed over a fixed period of time. If the true density $g(\cdot)$ belongs to the Nikol'ski class $\mathcal{N}(\beta,L)$, the estimator $\hat{g}_h(\cdot)$ reaches the rate of convergence $N^{(-2\beta/(2\beta+1))}$ under some conditions linking $T$ to $N$ and imposing that $T \rightarrow \infty$ when $N\rightarrow\infty$, such as $T \geq N^3$ in the multiplicative case \eqref{eq:mult.model.nonp} or $T \geq n^{5/2}$ in the linear case \eqref{eq:lin.model.nonp}. Unsurprisingly, the theoretical investigations on discretization and the numerical experiments reveal that the sampling interval needs to be small enough to preserve the performances of $\hat{g}_h(\cdot)$. In \cite{dion2016b}, the condition $n \Delta ^2 = o(1)$ on the sampling interval is explicitly stated. Let us stress that in contrast to the parametric methods previously exposed, tuning aspects such that the choice of the bandwidth value $h$ are at least as important as the sampling scheme to ensure good performance of the kernel estimator. The adjustment of the bandwidth has been theoretically investigated, but this is out of the scope of the present paper.  \\

Let us finish by mentioning that deconvolution methods built on estimators of the random effects have also been investigated for the nonparametric estimation of $g(\cdot)$. The idea has been introduced in Comte \& al. (2013) by noting that the estimators of the random effects could be of the form $\hat{\phi}_{i,T}=\phi_i + \xi_i$, with $(\xi_1,\ldots,\xi_N)$ \textit{i.i.d.} random variables representing some form of noise. This happens when the model includes a linear random effect in the drift (b), leading to estimator \eqref{eq:phiest2}, or a multiplicative random effect in the drift in case where $b(\cdot) =\sigma(\cdot)$. The approach was taken up and improved by Dion (2016) in a mixed-effects Ornstein-Uhlenbeck model,  which is a particular case of (b). The different forms for $\xi_i$ according to the model and the range of implemented cut-off strategies lead to different expressions of the density estimators, thus we do not give any formula. The significant advantage of the deconvolution approach is that it can be declined both in the case where $T$ is fixed and in the case where $T$ is large. Estimation rates similar to those of the kernel approach are obtained. As for the kernel approach, the sampling time interval plays a key role in the deconvolution performances. The impact of discretization on the deconvolution approach is discussed in both \cite{comte2013} and \cite{dion2016a} in addition to the order of the risk as $N$ is large.

\section{Summary and concluding remarks}
\label{sec:conclusion}

The present review paper addresses inference in mixed-effects SDE models. Although many estimation methods have been developed for these models, few contributions address the theoretical properties of the associated estimators. In this paper, we restrict our review to asymptotic results on estimation and we do not discuss algorithmic contributions. 

The most natural approach is to consider maximum likelihood estimation, but exact computation is rarely feasible in practice. As mentioned in Section \ref{sec:discretelik}, this is only possible under two conditions, first the knowledge of the transition densities of the SDEs in explicit form, and secondly the ability to integrate analytically the conditional likelihood of the processes on random effects. These two conditions are only met in very simple models; otherwise, other methods must be used to obtain estimators. Several variants of the maximum likelihood method have been proposed and theoretically studied: using a discretized version of the continuous-time likelihood, using the Euler approximation scheme of the processes, plugging estimators of the random effects in the random-effects likelihood. These methods lead to explicit contrasts for parameter estimation in linear SDEs with Gaussian and inverse Gamma random effects in the drift and in the diffusion coefficients respectively. It appears that except in the rare situations where these contrasts coincide with the exact likelihood of the processes, defining consistent and asymptotically Gaussian estimators necessarily requires coupling the asymptotics where the number of diffusions $N$ tends towards infinity with high-frequency observations of the individual trajectories.  Other distributions for the random effects have however not been studied since except for plug-in (see Section \ref{sec:plugin.parametric}) only specific choices would lead to explicit contrasts.  

To deal with non specific distributions for the random effects, nonparametric strategies have also been considered, but only in cases where the random effects entering the drift are one- or two-dimensional, and where the diffusion coefficient does not contain random effects. 

With the proposed estimation methods, either parametric or not, high-frequency observations are almost always necessary. Coupling diffusions with random effects lead to quite uncommon sampling schemes with regard to the literature, that sometimes relate the sampling step $\Delta$ with the number $N$ of processes ($N\Delta \rightarrow 0$ or $N\Delta^2 \rightarrow 0$ for instance) rather than to the number of observations per trajectory (such that $n\Delta^2 \rightarrow 0$). Collecting a large number of observations on a large number of trajectories is however not realistic in many applications, for instance in pharmacology (see discussions in \cite{donnet2013}), where a limited number of observations can be sampled in practice.

To our knowledge, asymptotic inference has exclusively been tackled in SDEs without measurement noise and with linear random effects in the drift and in the diffusion coefficient. Deriving asymptotic results while considering nonlinear random effects and/or noisy observations of the diffusions remain open questions. Some applications also use multidimensional diffusions with mixed-effects where not all coordinates of the processes are observed. To the best of our knowledge, there is no theory for partially observed SDEs with mixed-effects.

\section{Appendix: Observations and asymptotic frameworks for efficient parameter estimation in standard diffusion models}

Assume the following diffusion:
\begin{equation}
dX(t) = b(X(t),\varphi)dt + \sigma(X(t),\gamma)dW_t,
\label{eq:SDE}
\end{equation}
where $(W_t,t \geq 0)$ is a standard Brownian motion, $\varphi$ and $\gamma$ are unknown parameters, and the drift function $x \mapsto b(x,\varphi,t)$ and the diffusion coefficient function $x \mapsto \sigma(x,\gamma,t)$ are known up to parameters $\varphi$ and $\gamma$. We make the necessary assumptions to ensure the existence and uniqueness of a strong solution to \eqref{eq:SDE}.

In this appendix, it is considered that the process $X(t)$ is observed without measurement noise. The possible estimators for the parameters in the drift $\varphi$ and in the diffusion coefficient $\gamma$, and the asymptotic frameworks under which their theoretical properties can be studied depend on the nature of the observations. 

\subsection{Continuous-time observations}

Consider first the case where the process $(X(t))$ is continuously observed on a time interval $[0,T]$. Then we need to consider that $T\rightarrow \infty$, that the diffusion $(X(t))$ is ergodic, and to assume that the diffusion coefficient $\sigma(x,\gamma)$ is known, meaning that the value of parameter $\gamma$ is also known. Then, it is possible to estimate consistently the parameters $\varphi$ in the drift with the estimator $\hat{\varphi}_T$ maximizing the continuous-time likelihood of the process. Moreover, $\sqrt{T}(\hat{\varphi}_T - \varphi_0)$ converges in distribution to a Gaussian distribution when the time horizon $T$ goes to infinity (see \cite{kutoyants2004}).

\subsection{Discrete-time observations}

Without loss of generality, let us consider the case of discrete observations $X(t_1)$, \ldots , $X(t_n)$ of process $(X(t))$ with regular sampling $\Delta$ on a time interval $[0,T]$, where $t_j=j\Delta$, $j=1,\ldots,n$ and $T=n\Delta$. The estimation essentially use the Euler scheme of the diffusion. There are several possible asymptotic frameworks that lead to different theoretical properties for the parameter estimates. 
\begin{enumerate}
\item If the observation time interval $[0,T]$ is fixed, \textit{i.e.} $T<\infty$, then the asymptotics is obtained by considering the number of observations $n$ going to infinity while the sampling interval $\Delta=\Delta_n=T/n$ goes to zero. In this case, only the parameters in the diffusion coefficient $\gamma$ can be estimated consistently and the rate of convergence of the estimator $\hat{\gamma}_{n}$ is $n^{1/2}$.
\item If $T \rightarrow \infty$, the number of observations $n$ is systematically assumed to go to infinity. Both parameters $\varphi$ and $\gamma$ can be estimated but their asymptotic properties are strongly related are strongly related to the scheme of observation, the most common being 
\begin{enumerate}
\item a large sample scheme where the sampling interval $\Delta$ is fixed while $n \rightarrow \infty$. In this case, $\sqrt{n}$-consistent estimators can be derived for both the parameters in the drift $\varphi$ and in the diffusion coefficient $\gamma$.
\item a high-frequency scheme where the sampling interval $\Delta_n \rightarrow 0$ while $n \rightarrow \infty$ in such a way that $T=n \Delta_n \rightarrow \infty$. Under the condition $n\Delta_n^2\rightarrow 0$, $\sqrt{n\Delta_n}$ consistent estimators for $\varphi$ are derived, whereas for parameter $\gamma$, $\sqrt{n}$-consistent estimators are obtained. 
\end{enumerate}  
\end{enumerate}
We refer the reader to \cite{kessler2012} for more precise results. 

\section*{Conflict of interest}
On behalf of all authors, the corresponding author states that there is no conflict of interest.

\bibliographystyle{spbasic}      


\end{document}